\documentclass[journal]{IEEEtran}
\usepackage{graphicx}
\usepackage{amssymb}
\usepackage[cmex10]{amsmath}
\usepackage{pgfplots, filecontents}
\usepgfplotslibrary{patchplots}
\usepackage[]{algorithm2e}
\usepackage{enumerate}
\usepackage{color}
\usepackage{epstopdf}
\graphicspath{{./fig/}}
\usepackage{cite}
\usepackage{bm}

\newcommand{\ve}{\mathbf}
\newcommand{\m}{\mathbf}

\newcommand{\veh}[1]{\hat{\mathbf{#1}}}

\newcommand{\Lagr}{\mathcal{L}}

\usepackage{amsthm}
\newtheorem*{Proposition}{Result}

\ifCLASSINFOpdf
\else
\fi
\begin{document}
%
\title{Constrained Best Linear Unbiased Estimation}

%
%
%

\author{Oliver~Lang,~\IEEEmembership{Member,~IEEE,}
        Mario~Huemer,~\IEEEmembership{Senior Member,~IEEE,}
        and~Markus~Steindl
\thanks{O. Lang and M. Huemer are with the Institute of Signal Processing, Johannes Kepler University Linz, Linz, Austria, e-mail: oliver.lang@jku.at; mario.huemer@jku.at.}
\thanks{M. Steindl is with the Institute for Algebra, Johannes Kepler University Linz, Linz, Austria.}}

\maketitle

\begin{abstract}
The least squares (LS) estimator and the best linear unbiased estimator (BLUE) are two well-studied approaches for the estimation of a deterministic but unknown parameter vector. In many applications it is known that the parameter vector fulfills some constraints, e.g., linear constraints. For such situations the constrained LS estimator, which is a simple extension of the LS estimator, can be employed. In this paper we derive the constrained version of the BLUE. It will turn out that the incorporation of the linear constraints into the derivation of the BLUE is not straight forward as for the constrained LS estimator, but the final expression for the constrained BLUE is closely related to that of the constrained LS estimator. 
\end{abstract}


%
\IEEEpeerreviewmaketitle

\section{Introduction}
\label{sec:intro}

We consider classical estimation with an underlying linear model
\begin{equation}
\ve{y}  = \m{H} \ve{x} + \ve{n}, \label{equ:DECONV071}
\end{equation}
which is frequently used in many areas of signal processing. Here, $\ve{y} \in \mathbb{C}^{N_\ve{y} }$ is the vector of measurements, $\ve{x} \in \mathbb{C}^{N_\ve{x} }$ is a deterministic but unknown parameter vector, $\m{H} \in \mathbb{C}^{N_\ve{y} \times N_\ve{x}}$ is the measurement matrix (with varying requirements throughout the paper), and $\ve{n} \in \mathbb{C}^{N_\ve{y}}$ is zero mean measurement noise with known positive definite covariance matrix $\m{C}_{\ve{n}\ve{n}}$. The probability density function (PDF) of $\ve{n}$ is otherwise arbitrary. For this linear model, linear classical estimators such as the least squares (LS) estimator or the best linear unbiased estimator (BLUE) \cite{Kay-Est.Theory, Best_linear_unbiased_estimator_algorithm_for_received_signal-strength_based_localization, Unterrieder2015, Huemer2011, Pereira2013, Du2014, Noroozi2015} are well known and widely employed. 

In many applications it is known that the parameter vector fulfills some constraints, e.g., the linear constraints
\begin{align}
\m{A} \ve{x} = \ve{b}, \label{equ:Add_constraints_001h}	
\end{align}
with full row rank $\m{A} \in \mathbb{C}^{N_\ve{b} \times N_\ve{x}}$, $\ve{b} \in \mathbb{C}^{N_\ve{b} }$, $N_\ve{b} < N_\ve{x}$. One of many possible examples where a parameter vector fulfills such linear constraints is when $\ve{x}$ describes the impulse response of a linear time-invariant system that is unable to transmit a DC signal. Then, the sum of all elements in $\ve{x}$ must be zero, which can be described by appropriately choosing $\m{A}$ and $\ve{b}$. Since the parameter vector is assumed to fulfill \eqref{equ:Add_constraints_001h}, we seek for an estimator whose estimates $\veh{x}$ fulfill 
\begin{align}
\m{A} \veh{x} = \ve{b}. \label{equ:Add_constraints_001threeh}	
\end{align}
A modification of the LS estimator that fulfills \eqref{equ:Add_constraints_001threeh} can be found in standard textbooks, e.g., \cite{Kay-Est.Theory, Scharf1991} and is termed \emph{constrained} LS estimator. For full column rank $\m{H}$ implying $N_\ve{y} \geq N_\ve{x}$ the constrained LS estimator is formally given by
\begin{align}
 \veh{x}_{\mathrm{CLS}} ={}& \left( \m{I} - \m{Q}^{-1} \m{A}^H \left(  \m{A} \m{Q}^{-1} \m{A}^H \right)^{-1} \m{A} \right)  \m{Q}^{-1}  \m{H}^H \ve{y} \nonumber \\
 & + \m{Q}^{-1} \m{A}^H \left(  \m{A}\m{Q}^{-1} \m{A}^H \right)^{-1}\ve{b}, \label{equ:LS_add_const_008}
 \end{align}
where $ \m{Q} = \m{H}^H \m{H}$. Possible applications of this estimator can be found in \cite{Heinz2001, Altunbasak1997}. The constrained LS estimator can easily be derived by minimizing the standard LS cost function under the constraints \eqref{equ:Add_constraints_001h}, c.f. \cite{Kay-Est.Theory}. The cost function as well as the constraints are functions of $\ve{x}$. The cost functions of the BLUE \cite{Kay-Est.Theory}, however, are functions of the rows of the estimator matrix, and the incorporation of the constraints \eqref{equ:Add_constraints_001h} as a function of $\ve{x}$ is not straight forward as in the LS case. In this work we derive the constrained BLUE by converting the constraints \eqref{equ:Add_constraints_001h} into constraints fitting to BLUE's cost functions with the help of the nullspace of $\m{A}$. Finally, the resulting estimator as a function of linearly independent basis vectors of the nullspace of $\m{A}$ is converted to a form that shows great similarities with the constrained LS in \eqref{equ:LS_add_const_008}.

Outlook: In Sec.~\ref{sec:cerivation}, the constrained BLUE is derived, interestingly with weaker prerequisites than required for the constrained LS estimator in \eqref{equ:LS_add_const_008}. A simplified notation of the constrained BLUE, however, with more strict prerequisites is discussed in Sec.~\ref{sec:discussion}. A simulation example demonstrating the performance gain achievable with the constrained BLUE is provided in Sec.~\ref{sec:example}.

Notation: Lower-case bold face variables ($\ve{a}$, $\ve{b}$,...) indicate vectors, and upper-case bold face variables ($\m{A}$, $\m{B}$,...) indicate matrices. We further use $\mathbb{C}$ to denote the set of complex numbers, $(\cdot)^T$ to denote transposition, $(\cdot)^H$ to denote conjugate transposition, $(\cdot)^*$ to denote conjugation, $\m{I}^{n\times n}$ to denote the identity matrix of size $n\times n$, and $\m{0}^{m\times n}$ to denote the zero matrix of size $m\times n$. If the dimensions are clear from context we simply write $\m{I}$ and $\m{0}$. $E[\cdot]$ denotes the expectation operator, where we use a subscript to denote the averaging PDF.

\section{Derivation}
\label{sec:cerivation}

In coherence with the constrained LS estimator in \eqref{equ:LS_add_const_008}, we assume the estimator to be affine and of the form 
\begin{equation}
	\veh{x} = \m{E}\ve{y} + \ve{f}.  \label{equ:extensionsBLUEb021}	
\end{equation}
As the estimator is actually affine the term 'linear' in the abbreviation 'BLUE' might be somewhat misleading. However, since also for other affine estimators the term 'linear' is usually used, we call the estimator constrained BLUE. The goal is now to find the estimator matrix $\m{E}\in \mathbb{C}^{N_{\ve{x}} \times N_{\ve{y}}}$ and the vector $\ve{f}\in \mathbb{C}^{N_{\ve{x}}}$. The constrained BLUE has to fulfill two types of constraints. The first one is the unbiased constraint
\begin{align}
	E_{\ve{y}}[\hat{\ve{x}}] ={}& E_{\ve{y}}[ \m{E} \ve{y} + \ve{f}] = E_{\ve{n}}[ \m{E} ( \m{H}\ve{x} + \ve{n} )  + \ve{f}] \\
={}&	\m{E} \m{H}\ve{x}  + \ve{f} \stackrel{!}{=} \ve{x}. \label{equ:extensionsBLUEb005c} 
\end{align}
By letting $\ve{e}_i^H$ be the $i^\text{th}$ row of $\m{E}$, $x_i$ be the $i^\text{th}$ element of $\ve{x}$, and $f_i$ be the $i^\text{th}$ element of $\ve{f}$, 
the unbiased constraint for $\ve{e}_i^H$ can be extracted from \eqref{equ:extensionsBLUEb005c} and is of the form
\begin{equation}
	\ve{e}_i^H \m{H}\ve{x}  + f_i \stackrel{!}{=} x_i. \label{equ:extensionsBLUE013saa}
\end{equation}
The second type of constraints are given by \eqref{equ:Add_constraints_001threeh}. 
For each $i \in \{1, 2,..., N_\ve{x}\}$ the variance of $\hat{x}_i$ serves as a cost function which is a function of $\ve{e}_i$ given as 
\begin{align}
J(\ve{e}_i) ={}&   E_\ve{y} \left[ \left( \hat{x}_i - E_\ve{y} \left[ \hat{x}_i \right] \right) \left( \hat{x}_i - E_\ve{y} \left[ \hat{x}_i \right] \right)^H \right]  \\
 ={}& E_\ve{y} \Big[ \left(\ve{e}_i^H \ve{y} + f_i - \ve{e}_i^H E_\ve{y} \left[\ve{y}\right] - f_i\right)  \nonumber \\
 & \times \left(\ve{e}_i^H \ve{y} + f_i - \ve{e}_i^H E_\ve{y} \left[\ve{y}\right] - f_i\right)^H \Big]  \\
 ={}& E_\ve{n} \Big[ \left(\ve{e}_i^H \left( \m{H}\ve{x} + \ve{n} \right)  - \ve{e}_i^H \m{H}\ve{x} \right) \nonumber \\
 & \times \left(\ve{e}_i^H \left( \m{H}\ve{x} + \ve{n} \right)  - \ve{e}_i^H \m{H}\ve{x} \right)^H \Big] \\
  ={}& E_\ve{n} \left[ \left( \ve{e}_i^H \ve{n} \right) \left( \ve{e}_i^H \ve{n} \right)^H \right]  = \ve{e}_i^H \m{C}_{\ve{n}\ve{n}} \ve{e}_i. \label{equ:extensionsBLUE011}
\end{align}
We note, that \eqref{equ:Add_constraints_001threeh} and \eqref{equ:extensionsBLUE013saa} represent constraints in $\veh{x}$, however, the $i^\text{th}$ cost function is a function of the vector $\ve{e}_i$, which is conflicting. We are therefore now converting the constraints \eqref{equ:Add_constraints_001threeh} and \eqref{equ:extensionsBLUE013saa} into constraints in $\ve{e}_i$ and $f_i$ for $(i = 1,\hdots, N_\ve{x})$.

We start with an analysis of $\m{A} \ve{x} = \ve{b}$ in \eqref{equ:Add_constraints_001h}. This linear system of equations has an infinite number of solutions that can be described as 
\begin{align}
\ve{x} = \ve{x}_p + \ve{x}_1 \alpha_1 + \ve{x}_2 \alpha_2 + \hdots + \ve{x}_{N_0} \alpha_{N_0}, \label{equ:constr_BLUE_001}	
\end{align} 
where the vectors $\ve{x}_i$, $i = 1, \hdots, N_0$, span the nullspace of $\m{A}$ such that $\m{A} \ve{x}_i = \m{0}^{N_\ve{b} \times 1}$, $N_0$ is the dimension of the nullspace of $\m{A}$ with $N_0 = N_\ve{x} - N_\ve{b}$, the scalar coefficients $\alpha_i$, $i = 1, \hdots, N_0$ are in general complex valued and arbitrary, and $\ve{x}_p$ is an arbitrary particular solution of $\m{A} \ve{x} = \ve{b}$, e.g., the least norm solution $\ve{x}_p =  \m{A}^H \left( \m{A} \m{A}^H \right)^{-1}\ve{b}$. However, the particular choice of $\ve{x}_p$ is not of importance in the following. Eq.~\eqref{equ:constr_BLUE_001} can be brought in the form 
\begin{align}
\ve{x} = \ve{x}_p + \m{N} \ve{\bm{\alpha}}, \label{equ:constr_BLUE_003}	
\end{align} 
where 
\begin{align}
\m{N} = \begin{bmatrix} \ve{x}_1 & \hdots & \ve{x}_{N_0} \end{bmatrix} \in \mathbb{C}^{N_\ve{x} \times N_0}, \hspace{5mm} \ve{\bm{\alpha}} = \begin{bmatrix} \alpha_1 \\ \vdots \\ \alpha_{N_0} \end{bmatrix} \in \mathbb{C}^{N_0}. \label{equ:constr_BLUE_004}	
\end{align}
With this notation we have $\m{A}\m{N} = \m{0}^{N_\ve{b} \times N_0}$. Inserting \eqref{equ:constr_BLUE_003} into \eqref{equ:extensionsBLUEb005c} results in
\begin{align}
	E_{\ve{y}}[\hat{\ve{x}}] ={}& 	\m{E} \m{H} \left( \ve{x}_p + \m{N} \ve{\bm{\alpha}} \right)  + \ve{f} \stackrel{!}{=}\ve{x}_p + \m{N} \ve{\bm{\alpha}} \label{equ:extensionsBLUEb005g} \\
	\Leftrightarrow{}& \left( \m{E} \m{H} \m{N} - \m{N} \right)\ve{\bm{\alpha}} + \left( \m{E} \m{H} - \m{I} \right)\ve{x}_p  + \ve{f} \stackrel{!}{=} \ve{0}. \label{equ:extensionsBLUEb005h}
\end{align}
To fulfill this equation for every possible vector $\ve{\bm{\alpha}}$, we deduce the following two constraints for $ \m{E}$ and $\ve{f}$:
\begin{align}
	 \m{E} \m{H} \m{N} ={}& \m{N} \label{equ:extensionsBLUEb005i} \\
	\ve{f} ={}& \left( \m{I}- \m{E} \m{H} \right)\ve{x}_p. \label{equ:extensionsBLUEb005j}
\end{align}
Let the $i^\text{th}$ row of $\m{N}$ be denoted as $\ve{n}_i^H$, then the constraint for $\ve{e}_i^H$ can be extracted from \eqref{equ:extensionsBLUEb005i} and leads to
$\ve{e}_i^H \m{H} \m{N} = \ve{n}_i^H$. We are now finally able to formulate the constrained optimization problem for $\ve{e}_i$:
\begin{align}
&\ve{e}_{\text{CB},i} = \text{arg} \hspace{2mm} \underset{\ve{e}_i}{\text{min}} \hspace{2mm}\ve{e}_{i}^H \m{C}_{\ve{n}\ve{n}} \ve{e}_{i}   \hspace{5mm} \text{s.t.} \hspace{5mm} \ve{e}_i^H \m{H} \m{N} = \ve{n}_i^H. \label{equ:WB_Real010qqlidue}
\end{align}
We solve this constrained optimization problem using the complex valued Lagrangian multiplier method \cite{The_Complex_Gradient_Operator_and_the_CR-Calculus} and Wirtinger's calculus for deriving the complex valued gradients \cite{Wirtinger1927}. The Lagrangian cost function for this problem is given by
\begin{align}
	\Lagr(\ve{e}_i) ={}& \ve{e}_i^H \m{C}_{\ve{n}\ve{n}} \ve{e}_i + \ve{\bm{\lambda}}^H\left( \m{N}^H \m{H}^H \ve{e}_i - \ve{n}_i \right)  \nonumber \\
	&+ \ve{\bm{\lambda}}^T\left( \m{N}^T \m{H}^T \ve{e}_i^* - \ve{n}_i^* \right). \label{equ:extensionsBLUE014}
\end{align}
The Wirtinger derivative with respect to $\ve{e}_i$ produces
\begin{equation}
	\frac{\partial \Lagr(\ve{e}_i)}{\partial \ve{e}_i} = \ve{e}_i^H \m{C}_{\ve{n}\ve{n}}  + \ve{\bm{\lambda}}^H \m{N}^H \m{H}^H.  \label{equ:extensionsBLUE015}
\end{equation}
Setting \eqref{equ:extensionsBLUE015} equal to zero results in
\begin{equation}
	\ve{e}_{\text{CB},i}^H = - \ve{\bm{\lambda}}^H \m{N}^H \m{H}^H \m{C}_{\ve{n}\ve{n}}^{-1}. \label{equ:extensionsBLUE016}
\end{equation}
Assuming full column rank of $\m{H} \m{N}$, which implies $N_\ve{y} \geq N_0$, and inserting \eqref{equ:extensionsBLUE016} into the constraint in \eqref{equ:WB_Real010qqlidue} produces
\begin{align}
	- \ve{\bm{\lambda}}^H  ={}& \ve{n}_i^H \left( \m{N}^H \m{H}^H \m{C}_{\ve{n}\ve{n}}^{-1} \m{H} \m{N} \right)^{-1}. \label{equ:extensionsBLUE018} 
\end{align}
Reinserting this result into the expression for $\ve{e}_i^H$ in \eqref{equ:extensionsBLUE016} yields
\begin{equation}
	\ve{e}_{\text{CB},i}^H = \ve{n}_i^H \left( \m{N}^H \m{H}^H \m{C}_{\ve{n}\ve{n}}^{-1} \m{H} \m{N} \right)^{-1} \m{N}^H \m{H}^H \m{C}_{\ve{n}\ve{n}}^{-1}. \label{equ:extensionsBLUE019}
\end{equation}
Since $\ve{n}_i^H$ is the only term in \eqref{equ:extensionsBLUE019} that depends on the index $i$, the expression for the estimator matrix is given by 
\begin{equation}
	\m{E}_{\text{CB}} = \m{N} \left( \m{N}^H \m{H}^H \m{C}_{\ve{n}\ve{n}}^{-1} \m{H} \m{N} \right)^{-1} \m{N}^H \m{H}^H \m{C}_{\ve{n}\ve{n}}^{-1}. \label{equ:extensionsBLUE021}
\end{equation}
In the following, we denote $\m{P} = \m{H}^H \m{C}_{\ve{n}\ve{n}}^{-1} \m{H}$. Inserting \eqref{equ:extensionsBLUEb005j} and \eqref{equ:extensionsBLUE021} into \eqref{equ:extensionsBLUEb021} finally leads to the constrained BLUE in the form of
\begin{align}
	\veh{x}_{\text{CB}} &= \m{E}_{\text{CB}}\ve{y} + \left( \m{I}- \m{E}_{\text{CB}} \m{H} \right)\ve{x}_p \label{equ:extensionsBLUEb021s}	\\
	&= \m{N} \left( \m{N}^H \m{P}\m{N}  \right)^{-1}  \m{N}^H \m{H}^H \m{C}_{\ve{n}\ve{n}}^{-1} \ve{y} \nonumber \\
	& + \left( \m{I}- \m{N} \left( \m{N}^H \m{P}\m{N}  \right)^{-1}  \m{N}^H \m{P} \right)\ve{x}_p \label{equ:extensionsBLUEb021t}	\\	
	&=\m{N}\left( \m{N}^H\m{P}\m{N} \right)^{-1}  \m{N}^H \m{H}^H \m{C}_{\ve{n}\ve{n}}^{-1}  \left(  \ve{y} - \m{H}\ve{x}_p \right)  +  \ve{x}_p. \label{equ:constr_BLUE_031}
\end{align}
The covariance matrix of the constrained BLUE in \eqref{equ:constr_BLUE_031} can easily shown to be
\begin{align}
 \m{C}_{\veh{x}\veh{x},\mathrm{CB}} ={}&  \m{N}\left( \m{N}^H\m{P}\m{N} \right)^{-1}  \m{N}^H. \label{equ:easdxtensionsBLUEb021t}
 \end{align}

$\veh{x}_{\text{CB}}$ in \eqref{equ:constr_BLUE_031} is actually independent of the particular choice of $\ve{x}_p$. To prove this we first show that the identity 
 \begin{align}
 \m{T}  = \m{T} \m{A}^H \left( \m{A} \m{A}^H \right)^{-1} \m{A}, \label{equ:constr_BLUE_a032}
\end{align}
with 
\begin{align}
 \m{T} = \m{I}- \m{N} \left( \m{N}^H \m{P}\m{N}  \right)^{-1}  \m{N}^H \m{P},
\end{align}
holds. For that we utilize the matrix $[\m{A}^H \hspace{3mm} \m{N} ]$. Since $\m{A} \m{N} = \m{0}$, the column spaces of  $\m{A}^H$ and $\m{N}$ are orthogonal to each such that $[\m{A}^H \hspace{3mm} \m{N} ]$ is invertible. Multiplying \eqref{equ:constr_BLUE_a032} with $[\m{A}^H \hspace{3mm} \m{N} ]$ from the right results in $[ \m{T}\m{A}^H  \hspace{3mm} \m{0} ] = [\m{T}\m{A}^H \hspace{3mm} \m{0}  ] $. Since this equation is true and $[\m{A}^H \hspace{3mm} \m{N} ]$ is invertible, \eqref{equ:constr_BLUE_a032} is also true. Now replacing $\m{T} = \m{I}- \m{N} \left( \m{N}^H \m{P}\m{N}  \right)^{-1}  \m{N}^H \m{P}$ in the second line of \eqref{equ:extensionsBLUEb021t} by the right hand side of \eqref{equ:constr_BLUE_a032} gives
\begin{align}
	\veh{x}_{\text{CB}} &=  \m{N} \left( \m{N}^H \m{P}\m{N}  \right)^{-1}  \m{N}^H \m{H}^H \m{C}_{\ve{n}\ve{n}}^{-1} \ve{y} \nonumber \\
	& + \m{T} \m{A}^H \left( \m{A} \m{A}^H \right)^{-1} \m{A} \ve{x}_p \\
	&=  \m{N} \left( \m{N}^H \m{P}\m{N}  \right)^{-1}  \m{N}^H \m{H}^H \m{C}_{\ve{n}\ve{n}}^{-1} \ve{y} \nonumber \\
	& + \m{T} \m{A}^H \left( \m{A} \m{A}^H \right)^{-1} \ve{b}. \label{equ:extensionsBLUEb021tztr}	
\end{align}
That finally means that using any particular $ \ve{x}_p$ in \eqref{equ:extensionsBLUEb021t} yields the same result as using the least norm solution $\ve{x}_p =  \m{A}^H \left( \m{A} \m{A}^H \right)^{-1} \ve{b}$.

Another important note is that $N_\ve{y} \geq N_\ve{x}$ is not required for the application of \eqref{equ:constr_BLUE_031}, which is in contrast to the constrained LS estimator in \eqref{equ:LS_add_const_008}. In fact, the constrained BLUE in \eqref{equ:constr_BLUE_031} only requires full column rank of $\m{H}\m{N}$ in order for $\m{N}^H \m{P}\m{N} $ to be invertible. This implies that $N_\ve{y} \geq N_0$, but $N_\ve{y} < N_\ve{x}$ is allowed.

For the case that $N_\ve{y} \geq N_\ve{x}$ , $\m{H}$ has full column rank, and $\m{C}_{\ve{n}\ve{n}}$ is positive definite (as originally assumed) which implies that $\m{P}$ is invertible and positive definite, the expression for the constrained BLUE in \eqref{equ:constr_BLUE_031} can be simplified.

\section{Simplification and Discussion}
\label{sec:discussion}

Note that the expression of the constrained BLUE in \eqref{equ:constr_BLUE_031} requires the calculation of a basis of the nullspace of the matrix $\m{A}$. We will now derive an expression of the constrained BLUE that does not require this nullspace evaluation, but which requires $N_\ve{y} \geq N_\ve{x}$. With the assumptions of full column rank $\m{H}$ and positive definite $\m{C}_{\ve{n}\ve{n}}$ (as originally assumed) $\m{P}$ is invertible and positive definite, and the following identity holds:
\begin{align}
\m{N} & \left( \m{N}^H\m{P}  \m{N} \right)^{-1} \m{N}^H =  \nonumber \\
& \m{P}^{-1} - \m{P}^{-1}\m{A}^H \left( \m{A} \m{P}^{-1} \m{A}^H \right)^{-1}  \m{A} \m{P}^{-1}. \label{equ:constr_BLUE_032}
\end{align}
This identity can be proven the following way. The $i^\text{th}$ column of $\m{N}$ is denoted as $\ve{x}_i$ according to \eqref{equ:constr_BLUE_004}. Furthermore, the $i^\text{th}$ column of $\m{A}^H$ is denoted as $\ve{a}_i$. We first show that the vectors $\m{P}^{-1}\ve{a}_1,  \hdots, \m{P}^{-1}\ve{a}_{N_\ve{b}},  \ve{x}_1, \hdots, \ve{x}_{N_0} $ are linearly independent: Fix $c_1, \hdots, c_{N_\ve{b}}, d_i, \hdots, d_{N_0} \in \mathbb{C}$ such that
\begin{equation}
\sum_{i=1}^{N_\ve{b}} c_i \m{P}^{-1}\ve{a}_i + \sum_{j=1}^{N_0} d_i \ve{x}_i = \ve{0}. \label{equ:proof_001}
\end{equation}
For $\ve{u} = \sum_{i=1}^{N_\ve{b}} c_i \ve{a}_i $ and $\ve{v} = \sum_{j=1}^{N_0} d_i \ve{x}_i$ we have $\m{P}^{-1}\ve{u} + \ve{v} = \ve{0}^{N_{\ve{x}} \times 1}$. Left multiplication by  $\ve{u}^H$ yields $\ve{u}^H\m{P}^{-1}\ve{u} =0$ since $\ve{u}$ and $\ve{v}$ are orthogonal. Since $\m{P}^{-1}$ is positive definite, we have that $\ve{u}=\ve{0}$. By the linearly independence of all $\ve{a}_i$, all $c_i$ are $0$. By \eqref{equ:proof_001}, all $d_j$ are $0$. Thus the only solution of \eqref{equ:proof_001} is $c_i = d_j = 0$ for all $i, j$, or in other words $\m{P}^{-1}\ve{a}_1,  \hdots, \m{P}^{-1}\ve{a}_{N_\ve{b}},  \ve{x}_1, \hdots, \ve{x}_{N_0} $ are linearly independent. Hence, the square matrix $[ \m{P}^{-1}\m{A}^H  \hspace{3mm} \m{N} ]$ is invertible. Furthermore, the matrix $\m{B} = [\m{A}^H  \hspace{3mm}  \m{P}\m{N} ]$ is invertible. Right multiplying \eqref{equ:constr_BLUE_032} by $\m{B}$ yields $[\m{0}  \hspace{3mm} \m{N}] = [\m{P}^{-1}\m{A}^H \hspace{3mm} \m{N} ] - [\m{P}^{-1}\m{A}^H  \hspace{3mm} \m{0} ]$. Since this equation is true and $\m{B}$ is invertible, \eqref{equ:constr_BLUE_032} is also true.

Inserting  \eqref{equ:constr_BLUE_032} into \eqref{equ:constr_BLUE_031} finally yields
\begin{align}
 \veh{x}_{\mathrm{CB}} ={}& \left( \m{I} - \m{P}^{-1} \m{A}^H\left(\m{A} \m{P}^{-1} \m{A}^H \right)^{-1} \m{A} \right)  \m{P}^{-1} \m{H}^H  \m{C}_{\ve{n}\ve{n}}^{-1}\ve{y}  \nonumber \\
  &+ \m{P}^{-1} \m{A}^H \left(\m{A} \m{P}^{-1} \m{A}^H \right)^{-1}\ve{b}. \label{equ:constr_BLUE_019}
 \end{align}
For the constrained BLUE in \eqref{equ:constr_BLUE_019} one can easily show that the covariance matrix is
\begin{align}
 \m{C}_{\veh{x}\veh{x},\mathrm{CB}} ={}&  \m{P}^{-1} - \m{P}^{-1} \m{A}^H\left(\m{A} \m{P}^{-1} \m{A}^H \right)^{-1} \m{A} \m{P}^{-1}. \label{equ:constr_BLUE_019asd}
 \end{align}
The expression for the constrained BLUE in \eqref{equ:constr_BLUE_019} has the advantage that the nullspace of $\m{A}$ is not required. Furthermore, comparing the constrained LS estimator in \eqref{equ:LS_add_const_008} with the constrained BLUE in \eqref{equ:constr_BLUE_019} reveals that they are connected in a very similar way as it is the case for the LS estimator and the BLUE \cite{Kay-Est.Theory}. Finally, we end up with the following: 
\begin{Proposition} \label{prop:Constr_BLUE} Consider the linear model $\ve{y}  = \m{H} \ve{x} + \ve{n}$, where $\ve{y}\in \mathbb{C}^{N_\ve{y} }$ is the measurement vector, $\m{H}\in \mathbb{C}^{N_\ve{y} \times N_\ve{x}}$ is a known measurement matrix with $N_\ve{y} \geq N_\ve{x}$ and full column rank, and $\ve{n}\in \mathbb{C}^{N_\ve{y}}$ is a zero mean random noise vector with known positive definite covariance matrix $\m{C}_{\ve{n}\ve{n}}$. If $\ve{x}$ fulfills the linear constraints $\m{A} \ve{x} = \ve{b}$ with full row rank  $\m{A} \in \mathbb{C}^{N_\ve{b} \times N_\ve{x}}$, $\ve{b} \in \mathbb{C}^{N_\ve{b} }$, $N_\ve{b} < N_\ve{x}$, then the constrained BLUE minimizing the variances of the elements of $\hat{\ve{x}}_\mathrm{CB}$ such that $\hat{\ve{x}}_\mathrm{CB}$ fulfills $\m{A} \veh{x}_\mathrm{CB} = \ve{b}$ is given by \eqref{equ:constr_BLUE_019}. Its covariance matrix $\m{C}_{\hat{\ve{x}}\hat{\ve{x}},\mathrm{CB}}$ is given by \eqref{equ:constr_BLUE_019asd}. 

If $N_\ve{y} \geq N_\ve{x}$  does not hold, then let $\m{N}\in \mathbb{C}^{N_\ve{x} \times N_0}$ be the matrix built by linearly independent (column) basisvectors that span the nullspace of $\m{A}$. If $\m{H}\m{N}$ has full column rank (implying $Ny \geq N_0$), then the constrained BLUE for $\ve{x}$ fulfilling $\m{A} \veh{x}_\mathrm{CB} = \ve{b}$ is given by \eqref{equ:constr_BLUE_031}. Its covariance matrix $\m{C}_{\hat{\ve{x}}\hat{\ve{x}},\mathrm{CB}}$ is given by \eqref{equ:easdxtensionsBLUEb021t}.  
\end{Proposition}

\section{Example}
\label{sec:example}

We assume $\ve{x}\in \mathbb{C}^{5}$ to be the discrete-time impulse response of an unknown system. Additionally, we know that the system is unable to transmit any DC signals. Hence, the sum of all elements of $\ve{x}$ must be zero. This can be described by a linear constraint as in \eqref{equ:Add_constraints_001h}, where $\m{A}$ is a row vector of length $5$ with all elements being $1$, and where $\ve{b}$ is $0$. This example has the advantage that the constrained LS estimator and the constrained BLUE not only can be compared with the LS estimator and the BLUE but also with intuitive estimators as it will be demonstrated soon.

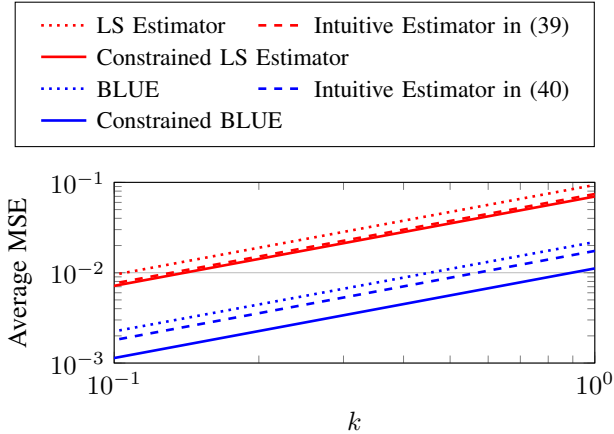
\begin{figure}[t!]
\begin{center}
\begin{tikzpicture}
\begin{loglogaxis}[compat=newest, 
width=0.9\columnwidth, height = .45\columnwidth, xlabel=$k$, 
ylabel style={align=center}, 
ylabel style={text width=3.4cm},
ylabel={Average MSE}, 
legend pos=north east, 
legend cell align=left,
legend columns=2, 
        legend style={
            /tikz/column 2/.style={
                column sep=5pt,
            },
        font=\small},
xmin = 0.1,
xmax = 1,
ymax = 0.1,
ymin = 10^(-3),
grid=major,
legend style={
at={(-0.25,1.3)},
anchor=north west}
]

\addplot[line width=1pt, color=red, style=dotted] table[x index =0, y index =4] {Average_MSE_x.dat};
\label{pc1}

\addplot[line width=1pt, color=red, style=solid] table[x index =0, y index =6] {Average_MSE_x.dat};
\label{pc3}

\addplot[line width=1pt, color=red, style=dashed] table[x index =0, y index =5] {Average_MSE_x.dat};
\label{pc2}

\addplot[line width=1pt, color=blue, style=dotted] table[x index =0, y index =1] {Average_MSE_x.dat};
\label{pc4}

\addplot[line width=1pt, color=blue, style=solid] table[x index =0, y index =3] {Average_MSE_x.dat};
\label{pc6}

\addplot[line width=1pt, color=blue, style=dashed] table[x index =0, y index =2] {Average_MSE_x.dat};
\label{pc5}

\end{loglogaxis}
\node [draw,fill=white,anchor=north east] at (rel axis cs: 2.00, 3.5) {\shortstack[l]{
	\begin{tabular}{ll}
		\ref{pc1} \small{LS Estimator} & \ref{pc2} \small{Intuitive Estimator in \eqref{equ:Clest_ex2_002}} \\
		\multicolumn{2}{l}{\ref{pc3} \small{Constrained LS Estimator}} \\
  		\ref{pc4} \small{BLUE}  &
  		\ref{pc5} \small{Intuitive Estimator in \eqref{equ:Clest_ex2_003}} \\
  		\multicolumn{2}{l}{\ref{pc6} \small{Constrained BLUE}} \\
\end{tabular}}} ;

\end{tikzpicture}
\caption{Average MSEs of the estimated impulse responses for various estimators. \label{fig:Clest_ex2_fig1}}
\end{center}
\end{figure}

The measurement vector $\ve{y}\in \mathbb{C}^{10}$ shall contain noisy measurements of the input samples convolved with the impulse response. The input samples written in vector form are denoted as $\ve{u}\in \mathbb{C}^{6}$. Thus, $\m{H}\in \mathbb{C}^{10 \times 5}$ is a convolution matrix built from the vector $\ve{u}$. The input samples are randomly drawn for every simulation run from a standard proper Gaussian distribution \cite{Complex-Valued-Signal-Processing-The-Proper-Way-to-Deal-With-Impropriety}. The covariance matrix of the complex proper noise vector $\ve{n}$ in \eqref{equ:DECONV071} is chosen as $\m{C}_{\ve{n}\ve{n}} = k \, \mathrm{diag}\{ [
1 \ \  1 \  \ 0.5 \  \ 0.5 \  \ 0.1 \  \ 0.1 \  \ 0.01 \  \ 0.01 \  \ 10^{-3} \  \ 10^{-3}
] \}$,\\ where $k$ is a scaling factor varied between $10^{-1}$ and $1$. The following estimators are considered: 
\begin{enumerate}
\item The LS estimator \cite{Kay-Est.Theory}, denoted as $\hat{\ve{x}}_{\mathrm{LS}}$,
\item the intuitive estimator resulting from subtracting the mean value from the estimates of the LS estimator 
\begin{equation}
\hat{\ve{x}} = \hat{\ve{x}}_{\mathrm{LS}} - \mathrm{mean}(\hat{\ve{x}}_{\mathrm{LS}}) \ve{1}^{N_\ve{x} \times 1}, \label{equ:Clest_ex2_002}
\end{equation}
where $ \ve{1}^{N_\ve{x} \times 1}$ denotes a column vector of length $N_\ve{x}$ with all elements being $1$,
\item the constrained LS estimator in \eqref{equ:LS_add_const_008},
\item the BLUE \cite{Kay-Est.Theory}, denoted as $\hat{\ve{x}}_{\mathrm{B}}$,
\item the intuitive estimator resulting from subtracting the mean value from the estimates of the BLUE
\begin{equation}
\hat{\ve{x}} = \hat{\ve{x}}_{\mathrm{B}} - \mathrm{mean}(\hat{\ve{x}}_{\mathrm{B}}) \ve{1}^{N_\ve{x} \times 1}, \label{equ:Clest_ex2_003}
\end{equation}
\item the constrained BLUE in \eqref{equ:constr_BLUE_019}.
\end{enumerate}
The resulting average MSEs (averaged over the elements of $\ve{x}$) plotted over $k$ are presented in Fig.~\ref{fig:Clest_ex2_fig1}. The LS estimator performs worst for all values of $k$. The estimation accuracy can be significantly increased by using the intuitive estimator in \eqref{equ:Clest_ex2_002}. An even better estimation accuracy is achieved by the constrained LS estimator. Even larger performance gains are obtained for the constrained BLUE when compared to the BLUE and the intuitive estimator in \eqref{equ:Clest_ex2_003}. Hence, compared to the LS-based estimators, it is even more beneficial to favor the constrained BLUE over the BLUE in this example.


\section{Conclusions}
\label{sec:concl}

This work closes a long existing gap in classical estimation theory, namely the derivation of the constrained BLUE for the case when the parameter vector fulfills linear constraints. We derived two versions of the constrained BLUE, the first one under even weaker prerequisites than for the well known constrained LS estimator, and the second one under similar prerequisites as the constrained LS estimator. The second version is also similar in form to the constrained LS solution.

%



\bibliography{References}

\end{document}